\documentclass[12pt,reqno]{article}
\usepackage{amssymb,amscd,amsmath,amsthm,color}
\newtheorem{theorem}{Theorem}[section]

\newtheorem{cor}[theorem]{Corollary}
\newtheorem{prop}[theorem]{Proposition}

\usepackage{enumerate,amssymb}
\theoremstyle{definition}

\newtheorem{question}[theorem]{Question}
\theoremstyle{remark}
\newtheorem{remark}[theorem]{Remark}

\numberwithin{equation}{section}

\def\bA{\mathbb{A}}
\def\bC{\mathbb{C}}

\def\bM{\mathbb{M}}

\def\bB{\mathbb{B}}

\begin{document}
\baselineskip=15pt

\title{A Pythagorean theorem for partitioned matrices}

\author{ Jean-Christophe Bourin{\footnote{Funded by the ANR Projet (No.\ ANR-19-CE40-0002) and by the French Investissements
 d'Avenir program, project ISITE-BFC (contract ANR-15-IDEX-03).}} \,and Eun-Young Lee{\footnote{This research was supported by
Basic Science Research Program through the National Research
Foundation of Korea (NRF) funded by the Ministry of
Education (NRF-2018R1D1A3B07043682)}  }   }

\date{ }

\maketitle

\vskip 10pt\noindent
{\small 
{\bf Abstract.} We establish a Pythagorean theorem for the absolute values of the blocks of a partitioned matrix. This leads to a series of  remarkable  operator  inequalities. For instance, if the  matrix $\bA$ is partitioned into three  blocks $A,B,C$, then
$$
|\bA|^3 \ge U|A|^3U^* + V|B|^3V^*+ W|C|^3W^*,
$$
$$
\sqrt{3} |\bA| \ge  U|A|U^* + V|B|V^*+ W|C|W^*,
$$
for some isometries $U,V,W$, and
$$
\mu_4^2(\bA) \le \mu_3^2(A) +\mu_2^2(B) + \mu_1^2(C)
$$
where $\mu_j$ stands for the $j$-th singular value.  Our theorem may be used to extend a result by Bhatia and Kittaneh for the Schatten $p$-noms and  to give a singular value version of Cauchy's Interlacing Theorem.


}

\section{Introduction}

Let $\bM_{d}$ denote the space of $d$-by-$d$ matrices. If $\bA\in\bM_{d}$, the polar decomposition holds, 
\begin{equation}\label{polar}
\bA = U |\bA|
\end{equation}
where $|\bA|\in\bM_{d}$ is positive semi-definite and $U\in\bM_{d}$ is a unitary matrix. The matrix $|\bA|$ is called the  absolute value of $\bA$, and its
  eigenvalues  are the singular values of $\bA$. The absolute value can be defined for $d\times d'$ matrices  $\bA\in\bM_{d,d'}$ as a positive matrix $|\bA|\in\bM_{d'}$, and the factor $U$ in \eqref{polar} is an isometry ($d\ge d'$) or a coisometry ($d<d'$).

 If $\bA$ is partitioned in some number of rectangular blocks, say four blocks $A,B,C,D$, it is of interest to have a relation between the absolute value $|\bA|$ and the absolute values of the blocks. By using the standard  inner product of $\bM_{d,d'}$, we immediately have the trace relation
$$
{\mathrm{Tr\,}} |\bA|^2 = {\mathrm{Tr\,}} |A|^2 + {\mathrm{Tr\,}} |B|^2+ {\mathrm{Tr\,}} |C|^2 +{\mathrm{Tr\,}} |D|^2.
$$
This note aims to point out a much stronger Pythagorean theorem, Theorem \ref{th-pyth}, and several consequences. This result holds for many partitionings of $\bA$, not only when $A$ is a block matrix partitioned into $p\times q$ blocks. For instance, one may consider the  matrix
$$
\bA=\begin{pmatrix} a_1&a_2 &b_1 &b_2 &b_3  \\  a_3&a_4 &b_4 &b_5 &b_6 \\  a_5&a_6 &c_1 &c_2 &d_1 \\ a_7&a_8 &c_3 &c_4 &d_2 \\ a_9&a_{10} &c_5 &c_6 &d_3
\end{pmatrix}
$$
partitioned into four obvious blocks $A,B,C,D$.

If $\bA$ is partitioned into $r$ blocks $A_k\in\bM_{n_k,m_k}$, we write
\begin{equation}\label{part}
\bA=\bigcup_{k=1}^r A_k =A_1 \cup \cdots \cup A_r
\end{equation}
where we can use the $=$ sign if one considers  $A_k$ not only as an element of $\bM_{n_k,m_k}$ but also as a submatrix of
$\bA$ with its position in $\bA$. 

We say that the partitioning \eqref{part}  is colum compatible, or that $\bA$ is partitioned into colum compatible blocks if for all pairs of indexes $k,l$, either $A_k$ and $A_l$ lie on the same set of columns of $\bA$, or $A_k$ and $A_l$ lie on two disjoint  sets of columns of $\bA$.  Similarly, \eqref{part}  is row compatible, if for all pairs of indexes $k,l$, either $A_k$ and $A_l$ lie on the same set of rows of $\bA$, or$A_k$ and $A_l$ lie on two disjoint  sets of rows of $\bA$. 

Our Pythagorean Theorem \ref{th-pyth} will be stated for row or column compatible blocks. An application is a Theorem of Bhatia and Kittaneh for the Schatten $p$-norms (Corollary \ref{corBK}). Another application is an inequality for the singular values of compression onto hyperplanes. A matrix $A\in\bM_d$ is an operator on $\bC^d$. Given a hyperplane ${\mathcal{S}}$  of $\bC^d$,  we have a unit vector $h$ such that ${\mathcal{S}}
=h^{\perp}$, that is $x\in{\mathcal{S}}\iff \langle h, x\rangle =h^*x=0$. The compression $ A_{\mathcal{S}}$ of $A$ onto ${\mathcal{S}}$ is the operator acting on ${\mathcal{S}}$ defined as the restriction of $EA$ to ${\mathcal{S}}$ where $E$ stands for the (orthogonal) projection onto ${\mathcal{S}}$. Theorem \ref{th-pyth} entails a bound for the singular values of $A_{\mathcal{S}}$ in terms of those of $A$. These results are given in Section 3; we state a special case in the following corollary. Let $\mu_j$ denote the $j$-th singular value arranged in nonincreasing order.

\vskip 5pt
\begin{cor}\label{normal} 
 Let $A\in\bM_d$ be a normal matrix and let  ${\mathcal{S}}$ be a hyperplane of  $\bC^d$ orthogonal to the unit vector $h$. Set $\beta= \| Ah\|^2 -|\langle h,Ah\rangle|^2$. Then, 
for $j=1,\ldots, d-1$,
$$
 \mu_j^2(A)\ge \mu_{j}^2(A_{\mathcal{S}}) \ge \mu^2_{j+1}(A) -\beta.
$$
\end{cor}

\vskip 5pt
We discuss the case of four and five blocks in Section 4.  For  four blocks, our Pythagorean theorem entails an interesting inequality stated in the next corollary.

\vskip 5pt
\begin{cor}\label{revtriangle}  Let $\bA\in\bM_{d,d'}$ be partitioned into  four  blocks $A, B, C, D$. Then, there exist some isometries $U,V,W,X$ of suitable sizes such that
$$
2|\bA| \ge U|A|U^* + V|B|V^* + W|C|W^* + X|D|X^*.
$$
\end{cor}

The last section is devoted to   several other operator inequalities such as the first inequality in the abstract.

\section{A Pythagorean theorem}

\vskip 5pt
\begin{theorem}\label{th-pyth}  Let $\bA\in\bM_{d,d'}$ be partitioned into  $r$ row or column compatible blocks $A_k\in\bM_{n_k, m_k}$. Then, there exist some isometries  $U_k\in\bM_{d',m_k}$ such that
$$
|\bA|^2= \sum_{k=1}^{r} U_k  |A_k|^2 U_k^*.
$$
\end{theorem}

\vskip 5pt Recall that  $U\in\bM_{d',m}$, $m\le d'$, is an isometry if $U^*U=\bold{1}_{m}$, the identity on $\bC^m$. If $\bA\in\bM_{d,1}$, then the theorem reads as Pythagoras Theorem.

\begin{proof} Consider a positive matrix in $\bM_{n+m}$ partitioned as
$$
\begin{bmatrix} A& X \\ X^* &B
\end{bmatrix}
$$
with diagonal blocks $A\in\bM_n$ and $B\in\bM_m$. By  \cite[Lemma 3.4]{BL} we have two unitary matrices $U,V\in\bM_{n+m}$ such that
\begin{equation}\label{key}
\begin{bmatrix} A& X \\ X^* &B
\end{bmatrix}
=U\begin{bmatrix} A& 0 \\ 0 &0
\end{bmatrix}U^* + V\begin{bmatrix} 0& 0 \\ 0 &B
\end{bmatrix}V^*,
\end{equation}
equivalently,
\begin{equation*}
\begin{bmatrix} A& X \\ X^* &B
\end{bmatrix}
=U_1AU_1^* +
U_2BU_2^*
\end{equation*}
for two isometry matrices $U_1\in\bM_{n+m,n}$ and $U_2\in\bM_{n+m,m}$. An obvious iteration of  \eqref{key} shows that, given a positive block matrix in $\bM_m$ partitioned into $p\times p$ blocks,
$$
\bB=\left(B_{i,j}\right)_{1\le i,j\le p},
$$
with square diagonal blocks $B_{i,i}\in \bM_{n_i}$ and $n_1+\cdots+n_p=m$, we have the decomposition
\begin{equation}\label{key2}
\bB=\sum_{i=1}^p U_i B_{i,i} U_i^*
\end{equation}
for some isometries $U_i\in\bM_{m, n_i}$.

We use \eqref{key2} to prove the theorem. Consider first the column compatible case. Thus we have a partitioning into $p$ block columns,
\begin{equation}\label{col}
\bA={\mathbf{C}}_1\cup\cdots\cup {\mathbf{C}}_p,
\end{equation}
and each block $A_k$ belongs to  one block column ${\mathbf{C}}_q$. By relabelling the $A_k$'s if necessary, we may assume that we have $p$ integers $1=\alpha_1<\alpha_2<\cdots<\alpha_p<r$ such that
$$
{\mathbf{C}}_q= A_{\alpha_q} \cup \cdots \cup A_{\alpha_{q+1}-1}, \quad 1\le q< p, \quad{\text{and}} \quad {\mathbf{C}}_p= A_{\alpha_p} \cup \cdots \cup A_{\alpha_{r}}.
$$
We also have a partitioning into $p$ block rows, 
\begin{equation}\label{row}
\bA^*={\mathbf{C}}^*_1\cup\cdots\cup {\mathbf{C}}^*_p,
\end{equation}
and multiplying \eqref{row} and \eqref{col} we then obtain a block matrix for $\bA^*\bA=|\bA|^2\in\bM_{d'}$,
$$
|\bA|^2= \left( {\mathbf{C}}^*_i{\mathbf{C}}_j\right)_{1\le, i,j\le p}.
$$
By the decomposition \eqref{key2} we have
$$
|\bA|^2=\sum_{i=1}^p U_i{\mathbf{C}}^*_i{\mathbf{C}}_i U_i^*
$$
for some isometries $U_i\in\bM_{d' , n_i}$, where $n_i$ is the number of columns of ${\mathbf{C}}_i$. Hence, with the convention $\alpha_{p+1}:=r+1$,
$$
|\bA|^2=\sum_{i=1}^p\sum_{k=\alpha_i}^{\alpha_{i+1}-1} U_i A_k^*A_k U_i^*
$$
establishing the theorem for a column compatible partitioning.

Now, we turn to the row compatible case. Thus we have a partitioning into $p$ block rows,
\begin{equation}\label{col1}
\bA={\mathbf{R}}_1\cup\cdots\cup {\mathbf{R}}_p,
\end{equation}
and each block $A_k$ belongs to  one block row ${\mathbf{R}}_q$ and, as in the colum compatible case, we may assume that we have $p$ integers $1=\alpha_1<\alpha_2<\cdots<\alpha_p<r$ such that
$$
{\mathbf{R}}_q= A_{\alpha_q} \cup \cdots \cup A_{\alpha_{q+1}-1}, \quad 1\le q< p, \quad{\text{and}} \quad {\mathbf{R}}_p= A_{\alpha_p} \cup \cdots \cup A_{\alpha_{r}}.
$$
We also have a partitioning into $p$ block columns, 
\begin{equation}\label{row1}
\bA^*={\mathbf{R}}^*_1\cup\cdots\cup {\mathbf{R}}^*_p
\end{equation}
Mutiply \eqref{row1} and \eqref{col1} and note that
\begin{equation}\label{note1}
|\bA|^2=\sum_{l=1}^p {\mathbf{R}}^*_l{\mathbf{R}}_l.
\end{equation}
with  $p$  block matrices in $\bM_{d'}$, ($l=1,\ldots,p$),
\begin{equation}\label{note2}
{\mathbf{R}}^*_l{\mathbf{R}}_l = \left( A^*_iA_j\right)_{\alpha_{l}\le i,j<\alpha_{l+1}}
\end{equation}
where we still use $\alpha_{p+1}:=r+1$. Applying the decomposition \eqref{key2} to the block matrices 
\eqref{note2} yields
$$
{\mathbf{R}}^*_l{\mathbf{R}}_l=\sum_{i=\alpha_l}^{\alpha_{l+1}-1} U_i |A_i|^2U_i^*
$$
for some isometries $U_i$ of suitable sizes, 
and combining with \eqref{note1} completes the proof. \end{proof}

Denote by $\mu_1(S) \ge \mu_2(S)\ge \cdots$ the singular values of a matrix $S\in\bM_{n,m}$. This list is often limited to $\min\{n,m\}$ elements, however we  can naturally define $\mu_k(S)=0$ for any index $k$ larger than $\min\{n,m\}$. Given two matrices of same size, a  classical inequality of Weyl asserts that
$$
\mu_{j+k+1} (S+T) \le \mu_{j+1}(S) +\mu_{k+1}(T)
$$
for all nonnegative integers $j$ and $k$. This inequality and Theorem \ref{th-pyth} entail the next corollary.

\vskip 5pt
\begin{cor}\label{cor-sing}  Let $\bA\in\bM_{d,d'}$ be partitioned into  $r$ row or column compatible blocks $A_k\in\bM_{n_k, m_k}$. Then, for all nonnegative integers $j_1,j_2,\ldots, j_r$,
$$
\mu^2_{j_1+j_2+\cdots+j_r+1}(\bA)\le \sum_{k=1}^{r} \mu^2_{j_k+1}(A_k).
$$
\end{cor}

\vskip 5pt
A special case of this inequality is given in the abstract for three blocks and $j_1=2$, $j_2=1$, $j_3=0$.

Since any partitinioning into three blocks is  row or column compatible we have the next corollary.

\vskip 5pt
\begin{cor}\label{cor-three} Let $\bA\in\bM_{d,d'}$ be partitioned into  three  blocks $A, B, C$. Then, there exist some isometries $U,V,W$ of suitable sizes such that
$$
|\bA|^2=U|A|^2U^* + V|B|^2V^* + W|C|^2W^*.
$$
\end{cor}

By using the triangle inequality for the Schatten $p$-norms we have the trace inequality
\begin{equation}\label{trace}
\left\{{\mathrm{ Tr}}\,|\bA|^{2p}\right\}^{1/p}\le \left\{{\mathrm{ Tr}}\,|A|^{2p} \right\}^{1/p}+ \left\{{\mathrm{ Tr}}\,|B|^{2p}\right\}^{1/p} + \left\{{\mathrm{ Tr}}\,|C|^{2p}\right\}^{1/p}, \quad p\ge 1,
\end{equation}
equivalently
\begin{equation}\label{schatt}
\|\bA\|^2_q \le \| A\|^2_q + \| B\|^2_q +\| C\|^2_q
\end{equation}
for all Schatten $q$-norms, $q\ge 2$.

Theorem \ref{th-pyth} entails another interesting relation between the blocks of a partitioned matrix and the full matrix.

\begin{cor}\label{cor-direct}  Let $\bA\in\bM_{d,d'}$ be partitioned into  $r$ row or column compatible blocks $A_k\in\bM_{n_k, m_k}$. 
Then, for some isometries $V_j\in\bM_{m,d'}$, with $m=\sum_{k=1}^r m_k$,
\begin{equation*}
\bigoplus_{k=1}^r |A_k|^2 =\frac{1}{r}\sum_{j=1}^r V_j|\bA|^2  V_j^*.
\end{equation*}
\end{cor}

\vskip 5pt
\begin{proof} From Theorem \ref{th-pyth} and the main result of \cite{BL-direct} we have
$$
\bigoplus_{k=1}^r U_k|A_k|^2U_k =\frac{1}{r}\sum_{j=1}^r W_j|\bA|^2  W_j^*
$$
for some isometries $U_k\in\bM_{d',m_k}$ and some isometries $W_j\in \bM_{rd',d'}$. Since
$$
\bigoplus_{k=1}^r |A_k|^2= C\left\{ \bigoplus_{k=1}^r U_k|A_k|^2U_k\right\}C^*
$$
for some contraction $C\in\bM_{m, rd'}$, we infer
$$
\bigoplus_{k=1}^r |A_k|^2=\frac{1}{r}\sum_{j=1}^r CW_j|\bA|^2  W_j^*C^*.
$$
If $|\bA|$ is invertible, then, taking trace, the above equality ensures that contractions $CW_j$ satisfy $W_j^*C^*CW_j={\mathbf{1}}_{d'}$ for all $j$. Hence the result is proved with $V_j=CW_j$. The general case follows by a limit argument.
\end{proof}

We are in a position to estimate the Schatten norms of the blocks with the full matrix.
The following corollary was first obtained by Bhatia and Kittaneh \cite{Bh-K} in case of a matrix partitioned into $n\times n$ blocks.

\vskip 5pt
\begin{cor}\label{corBK}  Let $\bA\in\bM_{d,d'}$ be partitioned into  $r$ row or column compatible blocks $A_k\in\bM_{n_k, m_k}$. 
Then, for all Schatten $q$-norms, $q\ge 2$,
\begin{equation*}r^{{\frac{2}{q}}-1}\sum_{k=1}^r \| A_k \|_q^2 \le \| \bA\|_q^2 \le \sum_{k=1}^r \| A_k \|_q^2
\end{equation*}
These two inequalities are reversed for $2>q>0$.
\end{cor}

\vskip 5pt
\begin{proof} For $p:=q/2\ge 1$, the second inequality contains \eqref{schatt} and immediately follows from  Theorem \ref{th-pyth} and the triangle inequality for the Schatten $p$-norms. Corollary \ref{cor-direct} gives
$$
\| |\bA|^2\|_p\ge \left\| |A_1|^2\oplus \cdots\oplus |A_r|^2 \right\|_p 
$$
and since the concavity of $t\mapsto t^{1/p}$ entails
$$
\left\| |A_1|^2\oplus \cdots\oplus |A_r|^2 \right\|_p = \left(\| |A_1|^2\|_p^p+\cdots  \| |A_r|^2\|_p^p\right)^{1/p}\ge r^{\frac{1}{p}-1}\left(\| |A_1|^2\|_p+\cdots+ \| |A_r|^2\|_p\right)
$$
we get the first inequality. These inequalities are reversed for $0<p<1$.
\end{proof}

\vskip 5pt
Corollary \ref{cor-direct} is relevant to Majorisation Theory. We take  this opportunity to point out an interesting
fact about majorisation in the next proposition. Though this result might be well-known of some experts, it does not seem to be in the literature. Let $A,B\in\bM_n^+$, the positive semi-definite cone of $\bM_n$. The majorisation $A\prec B$ means that
$$
\sum_{j=1}^k \mu_j(A) \le \sum_{j=1}^k \mu_j(B)
$$
for all $k=1,2,\ldots n$, with equality for $k=n$.
The majorisation $A\prec B$ is equivalent to
 $$A=\sum_{i=1}^{n} \alpha_iU_iBU_i^*$$
 for some unitary matrices $U_i\in\bM_n$ and weights $\alpha_i\ge 0$ with $\sum_{i=1}^{n} \alpha_i=1$. This can be easily derived from Caratheodory's theorem \cite{Zhan}. A more accurate statement holds.

\vskip 5pt
\begin{prop}\label{prop-maj}  Let $A,B\in\bM_n^+$, $A\prec B$.Then, for some unitary matrices $U_i\in\bM_n$,
$$
A=\frac{1}{n}\sum_{i=1}^n U_i BU_i^*.
$$
\end{prop}

\vskip 5pt
\begin{proof} By the Schur-Horn Theorem, we may assume that $A$ is the diagonal part of $B$. Then we use the simple idea of Equation (2) in the nice paper of Bhatia \cite{Bh-monthly}.
\end{proof}

\vskip 5pt
Note that Corollary \ref{cor-direct} may be restated as 
\begin{equation*}
\bigoplus_{k=1}^r |A_k|^2 =\frac{1}{r}\sum_{j=1}^r U_j(|\bA|^2\oplus O)  U_j^*.
\end{equation*}
for some unitary matrices $U_j$ and some fixed zero  matrix $O$. Hence we have an average of $r$ matrices in the unitary orbit of $|\bA|^2\oplus O$, this number $r$ being (much) smaller than the one given by Proposition \ref{prop-maj}, $d=m_1+\cdots +m_r$.

\section{Compression onto a hyperplane}

By a hyperplane of $\bC^d$ we mean a vector subspace of dimension $d-1$. The next corollary is a singular value version of Cauchy's Interlacing Theorem \cite[p.\ 59]{Bh}.

\vskip 5pt
\begin{cor}\label{hyper}  Let $A\in\bM_d$ and let  ${\mathcal{S}}$ be a hyperplane of  $\bC^d$ orthogonal to the unit vector $h$.
Set $\beta= \min\{\| Ah\|^2, \| A^*h\|^2 \} -|\langle h, Ah\rangle|^2$. Then, for all  $j=1,\ldots, d-1$,
$$
\mu_{j}^2(A) \ge \mu_{j}^2(A_{\mathcal{S}}) \ge \mu_{j+1}^2(A) -\beta.
$$
\end{cor}

\vskip 5pt
This double inequality is stronger than $\mu_{j}(A)\ge \mu_{j}(A_{\mathcal{S}}) \ge \mu_{j+1}(A) -\sqrt{\beta}$. If $A$ is a normal matrix, then $\| Ah\|=\| A^*h\|$ and we have Corollary \ref{normal}. If $A=V$ is a unitary matrix, $\mu_j(V)=1$ for all $j$ and $\| Vh\|=\| V^*h\|=1$ for all unit vectors, so we deduce from Corollary \ref{hyper} that 
$\mu_{j}(V_{\mathcal{S}}) \ge |\langle h, Vh\rangle|$. In fact one can easily check that $\mu_j(V_{\mathcal{S}})=1$  for $j\le d-2$ and $\mu_{d-1}(V_{\mathcal{S}})= |\langle h, Vh\rangle|$. Hence Corollary \ref{hyper} is sharp.

\vskip 5pt
\begin{proof} (Corollary \ref{hyper}) The inequality $\mu_j(A)\ge \mu_j(A_{\mathcal{S}})$ is trivial. To deal with the other inequality 
we may assume that $h$ is the last vector of the canonical basis and that  $A_{\mathcal{S}}$ is the submatrix of $A$ obtained by deleting the last column and the last line. We partition $A$ as
$$
A=A_{\mathcal{S}}\cup B\cup C
$$
where $B$ contains  the $d-1$ entries below $A_{\mathcal{S}}$ and $C$ is the last colum of $A$. We then apply to this partitioning 
Corollary \ref{cor-sing} with  $j_1=j-1$, $j_2=0$, and $j_3=1$ to get
$$
\mu_{j+1}^2(A) \le \mu_{j}^2(A_{\mathcal{S}}) +  \mu_1^2(B) +\mu_2^2(C).
$$
Since $\mu_2(C)=0$, we have
\begin{equation*}\
\mu_{j+1}^2(A) -\mu_1^2(B) \le \mu_{j}^2(A_{\mathcal{S}})
\end{equation*}
Observe that $\mu_1^2(B)= \| A^*h\|^2  -|\langle hAh\rangle|^2$, hence
\begin{equation}\label{svi1}
\mu_{j}^2(A_{\mathcal{S}}) \ge \mu_{j+1}^2(A) -\| A^*h\|^2  +|\langle h, Ah\rangle|^2.
\end{equation}

We may also partition $A$ as
$$
A=A_{\mathcal{S}}\cup R\cup L
$$
where $R$ contains the $d-1$ entries on the right of $A_{\mathcal{S}}$ and $L$ stands for the last line of $A$.
Arguing as above with $R$ and $L$ in place of $B$ and $C$ yields
\begin{equation}\label{svi2}
\mu_{j+1}^2(A) -\mu_{j}^2(A_{\mathcal{S}}) \le \mu_1^2(R)= \| Ah\|^2  -|\langle h, Ah\rangle|^2
\end{equation}
Combining \eqref{svi1} and \eqref{svi2} completes the proof.
\end{proof}

\vskip 5pt
\begin{cor}\label{cor-var}  Let $\bA\in\bM_{d,d'}$ be partitioned into  $r$ row or column compatible blocks $A_k\in\bM_{n_k, m_k}$. Then, for each block $A_k$ and all  $j\ge 1$,
$$
\mu^2_j(\bA)-\mu^2_j(A_k)\le \sum_{l\neq k} \mu^2_1(A_l).
$$
\end{cor}

\vskip 5pt
\begin{proof} Apply Corollary \ref{cor-sing} with $j_k=j-1$ and $j_l=0$ for all $l\neq k$. \end{proof}

\vskip 5pt
\begin{cor}\label{cor-comp}   Let $A\in\bM_d$ and let  ${\mathcal{S}}$ be a hyperplane of  $\bC^d$ orthogonal to the unit vector $h$. Then for all $j=1,\ldots, d-1$,
$$
\mu^2_j(A)-\mu^2_j(A_{\mathcal{S}})\le \| Ah\|^2 + \| A^*h\|^2-|\langle h,Ah\rangle|^2.
$$
\end{cor}

\vskip 5pt
\begin{proof}
We may suppose that $h$ is the last vector of the canonical basis and we  partition $A$ into three blocks : $A_{\mathcal{S}}$, the last column of $A$, and the $d-1$ entries below $A_{\mathcal{S}}$. We then apply the previous corollary.
 \end{proof}

\section{Four and five blocks}

Partitionings into four blocks are not necessarily row or column compatible. However, for such partitionings, Theorem \ref{th-pyth} still holds.

\vskip 5pt
\begin{cor}\label{cor-four} Let $\bA\in\bM_{d,d'}$ be partitioned into  four  blocks $A, B, C, D$. Then, there exist some isometries $U,V,W,X$ of suitable sizes such that
$$
|\bA|^2=U|A|^2U^* + V|B|^2V^* + W|C|^2W^* + X|D|^2X^*.
$$
\end{cor}

\vskip 5pt
\begin{proof} We assume that $A$ is the block in the upper left corner and we distinguish three cases.

(1) 
 $A$ has the same number  $d$ of  lines as $\bA$.
In such a case, letting $A'=B\cup C\cup D$, the partitioning $\bA=A\cup A'$ is  column compatible, and we have two isometry matrices $U, U'$ such that 
\begin{equation}\label{f1}
|\bA|^2=U|A|^2U^* + U'|A'|^2U'^*.
\end{equation}
Since $A'$ is partitioned into three blocks, necessarily a row or column partitioning,  we can apply the theorem to obtain the decomposition
\begin{equation}\label{f2}
|A'|^2=V'|B|^2V'^* + W'|B|^2W'^* +X'|B|^2X'^* 
\end{equation} 
for some isometry matrices $V',W',X'$ of suitable sizes. Combining \eqref{f1} and \eqref{f2} we get the conclusion of the corollary with the isometry matrices $V=U'V'$, $W=U'W'$, and $X=U'X'$.

(2) 
 $A$ has the same number  $d'$ of  columns as $\bA$.
Letting again $A'=B\cup C\cup D$, the partitioning $\bA=A\cup A'$ is  row compatible, and we may argue as in case (1).

(3)
$A$ has  $l<d$ lines and  $c<d'$ columns. There exist then a block, say $B$, on the top position, and just on the right of $A$, and another block, say $C$ just below $A$ and on the left side.   We consider three subcases (a), (b), (c).

(a) $B$ has less than $l$ lines. Then, the last block $D$ is necessarily below $B$ with the same number of columns as $B$, and so $C$ has the same number of columns as $A$, hence $\bA=A\cup C \cup B \cup D$ is a column compatible partitioning  and we can apply  the theorem.

(b) $B$ has exactly $l$ lines, like $A$. We denote by $\gamma$  the number of columns of $B$ and we consider three situations.

(i) $C$ has more than $c+\gamma$ columns. Then necessarily $C$ has $d'$ columns and $D$ is the upper right block with $l$ lines, hence $\bA=A\cup B \cup  D \cup C $ is  a line compatible partitioning  and we may apply the theorem.

(ii)  $C$ has exactly $c+\gamma$ columns. Then, letting $\bA''=A\cup B\cup C$ with have a partitioning into three blocks, and $\bA = \bA''\cup D$. Thus applying the theorem twice as in case (1) yields the conclusion.

(iii) $C$ has fewer than $c+\gamma$ columns. Then $D$ is the lower right block, with the same number of lines as $B$, and $\bA$ is partitioned into line compatible blocks. Thus the theorem can be applied.

(c)   $B$ has more lines than $A$.  Let $\lambda$ be the number of line of $B$. Hence $\lambda>l$. There exist two situations

(iv) $\lambda<d$. Then $D$ is the lower right block, with the same number of columns as $B$,  and $\bA$ is partitioned into line compatible blocks. Thus we may apply the theorem.

(v) $\lambda=d$. Then $A'''= A\cup C\cup D$ is a partitioning into three blocks and $\bA= A'''\cup B$, thus applying twice the theorem completes the proof.
\end{proof}

We do not know wether Corollary \ref{cor-four} can be extended or not to any partitioning in five blocks. For instance we are not able to prove or disprove a version of Corollary \ref{cor-four} for the matrix 
$$
\bA=\begin{pmatrix} a_1&a_2 &a_3 &b_1 &b_2  \\  a_4&a_5 &a_6 &b_3 &b_4 \\  d_1&d_2 &x &b_5 &b_6 \\ d_3&d_4 &c_1&c_2 &c_3 \\ d_5&d_6 &c_4&c_5 &c_6 
\end{pmatrix}
$$
partitioned into five obvious blocks $A,B,C,D,X$. Hence, that Theorem \ref{th-pyth} holds or not for  any partitioning into five blocks is an open problem. More generally, we may consider the following two questions.

\begin{question} For which partitionings does Theorem \ref{th-pyth} hold ? For which partitionings does Corollary \ref{cor-sing} hold ?
\end{question}

\vskip 5pt
Matrices partitioned into four blocks (usually of same size) are  comon examples of partitionings. A nontrivial inequality follows from the previous corollary.

\vskip 5pt
\begin{cor}\label{cor-four2} Let $\bA\in\bM_{d,d'}$ be partitioned into  four  blocks $A, B, C, D$, and let  $p>2$. Then, there exist some isometries $U,V,W,X$ of suitable sizes such that
$$
2^{2-p}|\bA|^{\color{red}p}\le U|A|^pU^* + V|B|^pV^* + W|C|^pW^* + X|D|^pX^*.
$$
The inequality reverses for  $2>p>0$.
\end{cor}

\vskip 5pt
Letting $p=1$ we have  Corollary \ref{revtriangle} with the constant 2 which is sharp, even for a positive block matrix,  as shown by the simple example 
$$
\bA=\begin{bmatrix} A&A \\ A&A\end{bmatrix}.
$$

\vskip 5pt
\begin{proof}
For any monotone convex function $f(t)$ on the nonnegative axis, we have thanks to \cite[Corollary 2.4]{BL} and Corollary  \ref{cor-four},
\begin{align*}
f\left(\frac{|\bA|^2}{4}\right) &=f\left(\frac{U|A|^2U^* + V|B|^2V^* + W|C|^2W^* + X|D|^2X^*}{4}\right)  \\
&\le \Lambda\frac{f(U|A|^2U^*)+ V(f|B|^2V^*) + f(W|C|^2W^*) + f(X|D|^2X^*)}{4}\Lambda^*
\end{align*}
for some unitary matrix $\Lambda\in\bM_d'$. Picking $f(t)=t^{p/2}$ with $p>2$ yields the result. The reverse inequalities hold for monotone concave functions $f(t)$  and $0<p<2$.
\end{proof}

\vskip 5pt
\begin{remark} The version of Corollary \ref{cor-four2} for three blocks $A,B,C$, and $p=1$ reads as the  inequality of the abstract,
$$
\sqrt{3} |\bA| \ge  U|A|U^* + V|B|V^*+ W|C|W^*.
$$
The constant $\sqrt{3}$ is the best one: we cannot take a smaller constant for
 $$\bA=\begin{bmatrix} x&y&z \\ x&y&z\\ x&y&z\end{bmatrix}$$
 partitioned into its three lines. For two blocks, a similar sharp inequality holds with the constant $\sqrt{2}$.
\end{remark}

\section{Concave or convex functions}

For sake of simplicity we state our results for a square matrix $\bA$ partitioned  into blocks. By adding some zero rows or zero columns to a rectangular matrix, we could obtain statements for rectangular matrices (Remark \ref{last}).

Suppose that $\bA\in\bM_d$ is partitioned into blocks $A_k\in\bM_{n_k,m_k}$, $k=1,\ldots, r$. From Thompson's triangle inequality
(\cite{T} or \cite[p.\ 74]{Bh} we have
\begin{equation}\label{thomp}
|\bA| \le \sum_{k=1}^r U_k|A_k|U_k^*
\end{equation}
for some isometry matrices $U_k\in\bM_{d,m_k}$. The equality of Theorem \ref{th-pyth} and \eqref{thomp} suggest several other  inequalities, in particular,  if $\bA$ is partioned in row or column compatible bloks,
\begin{equation}\label{spe}
|\bA|^3 \ge \sum_{k=1}^r V_k|A_k|^3V_k^*
\end{equation}
for some isometries $V_k\in\bM_{d,m_k}$. This is indeed true as shown in the following theorem. We do not know if \eqref{spe} can be extended to any partitioning. Corollary \ref{cor-four} and the proof of Theorem  \ref{th-convex} show that \eqref{spe} holds for four blocks. The case of five blocks is open.

\vskip 5pt
\begin{theorem}\label{th-convex}  Let $\bA\in\bM_{d}$ be partitioned into  $r$ row or column compatible blocks $A_k\in\bM_{n_k, m_k}$, and let $\psi(t)$ be a monotone function on $[0,\infty)$ such that $\psi(\sqrt{t})$ is convex and $\psi(0)= 0$. Then there exist some isometries $V_k\in\bM_{d,m_k}$ such that
$$
\psi(|\bA|)\ge \sum_{k=1}^r V_k  \psi(|A_k|) V_k^*.
$$
\end{theorem}

\vskip 5pt
Theorem \ref{th-convex} considerably improves \eqref{trace}. A special case with $\psi(t)=t^3$ is given in the abstract.

\vskip 5pt
\begin{proof} Let $g(t)$ be a monotone convex function on $[0,\infty)$ such that $g(0)\le 0$, and let $A,B\in\bM_n$ be positive (semidefinite). By \cite{AB} or \cite[Corollary 3.2]{BL} we have
\begin{equation*}
g(A+B) \ge Ug(A)U^* + Vg(B)V^*
\end{equation*}
for some unitary matrices $U,V\in\bM_n$. Using this inequality and Theorem \ref{th-pyth} we infer
\begin{equation*}\label{bien}
g(|\bA|^2)\ge \sum_{k=1}^r W_k g(U_k|A_k|^2U_k^*) W_k^*
\end{equation*}
for some unitary matrices $W_k$ and some isometry matrices $U_k\in\bM_{d,m_k}$. If $g(0)=0$, we have
$g(U_k|A_k|^2)U_k=U_kg(|A|^2)U_k^*$. Hence
\begin{equation*}
g(|\bA|^2)\ge \sum_{k=1}^r V_k g(|A_k|^2) V_k^*
\end{equation*}
with the isometry matrices $V_k=W_kU_k$. Applying this to $g(t)=\psi(\sqrt{t})$ completes the proof.
\end{proof}

\vskip 5pt
\begin{cor}\label{cor-concave}  Let $\bA\in\bM_{d}$ be partitioned into  $r$ row or column compatible blocks $A_k\in\bM_{n_k, m_k}$, and let $\varphi(t)$ be a nonnegative  function on $[0,\infty)$ such that $\varphi(\sqrt{t})$ is concave. Then there exist some isometries $U_k\in\bM_{d,n_k}$ such that
$$
\varphi(|\bA|)\le \sum_{k=1}^{r} U_k  \varphi(|A_k|) U_k^*.
$$
\end{cor}

\vskip 5pt
\begin{proof} Since $\varphi(\sqrt{t})$ is nonnegative and concave, it is necessarily a monotone function (nondecreasing), hence continuous on $(0,\infty)$. Since we are dealing with matrices we may further suppose that $\varphi(t)$ is also continuous at $t=0$.

(1) Assume that $\varphi(0)=0$. Theorem \ref{th-convex} applied to $\psi(t)=-\varphi(t)$ proves the corollary.

(2) Assume that $\varphi(0)>0$. Since the continuous functional calculus is continuous on the positive semidefinite cone of any $\bM_m$, by a limit argument, we may assume that $|\bA|$ is invertible. So, suppose that the spectrum of $|\bA|^2$ lies in an interval $[r^2, s^2]$ with $r > 0$.
Define a convex function $\phi(\sqrt{t})$ by $\phi(\sqrt{t})=\varphi(\sqrt{t})$ for $t\ge r^2$, $\phi(0)=0$, and the graph of  $\phi(\sqrt{t})$ on $[0,r^2]$ is a line segment. Hence $\phi(t)\le \varphi(t)$ and  $\phi(|\bA|)=\varphi(|\bA|)$. Applying case (1) to $\phi$ yields
$$
\varphi(|\bA|)=\phi(|\bA|)\le \sum_{k=1}^{r} U_k  \phi(|A_k|) U_k^* \le \sum_{k=1}^{r} U_k  \varphi(|A_k|) U_k^*
$$
for some isometry matrices $U_k$.
\end{proof}

\vskip 5pt
The next three corollaries follow from Corollary \ref{cor-concave}.

\vskip 5pt
\begin{cor}\label{cor-1}  Let  $\bA\in\bM_{mn}$ be partitioned into a family of $m\times m$ blocks $A_{i,j}\in\bM_{n, n}$, and let $0< q\le 2$. Then there exist some isometries $U_{i,j}\in\bM_{mn, n}$ such that
$$
|\bA|^q\le \sum_{i,j=1}^m U_{i,j}  |A_{i,j}|^q U_{i,j}^*.
$$
\end{cor}

\vskip 5pt
\begin{cor}\label{cor-2}  Let  $\bA\in\bM_{mn}$ be partitioned into a family of $m\times m$ blocks $A_{i,j}\in\bM_{n, n}$, let $s\ge 1$ and  $0< q\le 2$. Then,
$$
\left\{{\mathrm{Tr}\,}|\bA|^{qs}\right\}^{1/s}\le \sum_{i,j=1}^m  \left\{{\mathrm{Tr}\,}|A_{i,j}|^{qs}\right\}^{1/s}.
$$
\end{cor}

\vskip 5pt
\begin{cor}\label{cor-3}  Let  $A\in\bM_{n}$, let $c_k$ be the  norm of the $k$-th column of $A$ and let $0< q\le 2$. Then there exist some rank one projection $E_k\in\bM_{n}$ such that
$$
|\bA|^q\le \sum_{k=1}^n c_k^q E_k
$$
\end{cor}

\vskip 5pt
The last corollaries follow from Theorem \ref{th-convex}.

\vskip 5pt
\begin{cor}\label{cor-4}  Let  $\bA\in\bM_{mn}$ be partitioned into a family of $m\times m$ blocks $A_{i,j}\in\bM_{n, n}$, and let $p\ge 2$. Then there exist some isometries  $U_{i,j}\in\bM_{mn, n}$ such that
$$
|\bA|^p\ge \sum_{i,j=1}^m U_{i,j}  |A_{i,j}|^p U_{i,j}^*.
$$
\end{cor}

\vskip 5pt
\begin{cor}\label{cor-5}  Let  $\bA\in\bM_{mn}$ be partitioned into a family of $m\times m$ blocks $A_{i,j}\in\bM_{n, n}$, let $0\le s\le 1$ and  $p\ge 2$. Then,
$$
\left\{{\mathrm{Tr}\,}|\bA|^{ps}\right\}^{1/s}\ge \sum_{i,j=1}^m  \left\{{\mathrm{Tr}\,}|A_{i,j}|^{ps}\right\}^{1/s}.
$$
\end{cor}

\vskip 5pt
\begin{cor}\label{cor-6}  Let  $A\in\bM_{n}$, let $r_k$ be the  norm of the $k$-th row of $A$ and let $p\ge 2$. Then there exist some rank one projections $E_k\in\bM_{n}$ such that
$$
|\bA|^p\ge \sum_{k=1}^n c_k^p E_k.
$$
\end{cor}

\vskip 5pt
\begin{remark}\label{last} The proof of Theorem \ref{th-convex} is the same for a $d\times d'$ matrix $\bA$. So Corollary \ref{cor-concave} also holds for $\bA\in\bM_{d,d'}$ if $\varphi(0)=0$. In case of $d\ge d'$, we may again use a limit argument and assume that $|\bA|$ is invertible. In case of $d'>d$ we may argue as follows. Add some zero lines to $\bA$ in order to obtain
a square matrix $\bA_0\in\bM_d'$. Let $B_1\ldots,B_p$ be the blocks at the bottom of $\bA$, and $R_1,\ldots,R_q$ be the remaining blocks of $\bA$. Add some zeros to the blocks $B_i$ in order to obtain blocks $B_i^0$ of $\bA_0$ in such a way that
$$
\bA_0 =\left(\bigcup_{i} B_i^0\right)\cup \left(\bigcup_{j} R_j\right)
$$
is a row or colum compatible partitioning of $\bA_0$. Since it is a square matrix, we may apply Corolloray \ref{cor-concave} and since
$|\bA_0|=|\bA|$ and  $|B_i^0|=|B^i|$, we see that Corollary \ref{cor-concave} holds for $d\times d'$ matrices.
\end{remark}

\vskip 15pt
\noindent
Jean-Christophe Bourin

\noindent
Laboratoire de math\'ematiques, 

\noindent
Universit\'e de Bourgogne Franche-Comt\'e, 

\noindent
25 000 Besan\c{c}on, France.

\noindent
Email: jcbourin@univ-fcomte.fr

  \vskip 10pt
\noindent
Eun-Young Lee

\noindent Department of mathematics, KNU-Center for Nonlinear
Dynamics,

\noindent
Kyungpook National University,

\noindent
 Daegu 702-701, Korea.

\noindent Email: eylee89@knu.ac.kr


\begin{thebibliography}{99}
{\small

 
\bibitem{AB}  J.S.\ Aujla and J.-C.\ Bourin, Eigenvalue inequalities for convex and log-convex functions,
\textit{Linear Algebra Appl.} 424 (2007), 25--35.

 \bibitem{Bh} R.\ Bhatia, Matrix Analysis, Gradutate Texts in Mathematics, Springer, New-York, 1996.

 \bibitem{Bh-monthly}  R.\ Bhatia, Pinching, trimming, truncating, and averaging of matrices. {\it Amer.\ Math.\ Monthly}  107  (2000),  no. 7, 602--608. 

\bibitem{Bh-K} R.\ Bhatia and  F.\ Kittaneh, Norm inequalities for partitioned operators and an application. {\it  Math. Ann}.\  287  (1990),  no.\ 4, 719--726.

\bibitem{BL} J.-C.\ Bourin and E.-Y.\ Lee,  Unitary orbits of Hermitian operators with convex or concave functions,  {\it Bull.\ Lond.\ Math.\ Soc.}\ 44 (2012), no.\ 6, 1085--1102. 


\bibitem{BL-direct} J.-C.\ Bourin and E.-Y.\ Lee,  Direct sums of positive semi-definite matrices. {\it Linear Algebra Appl}.\  463  (2014), 273--281.

\bibitem{T} R.-C.\ Thompson, Convex and concave functions of singular values of matrix sums, {\it Pacific J.\ Math.}\ 66 (1976), 285--290.

\bibitem{Zhan} X.\ Zhan, The sharp Rado theorem for majorizations, {\it Amer.\ Math.\ Monthly} 110
(2003) 152--153.


}

\end{thebibliography}
\end{document}